\documentclass{ws-ijgmmp}

\usepackage{graphicx}
\usepackage{amsmath, amssymb, bm}
\usepackage{mathtools}
\usepackage{subfigure}
\usepackage{wrapfig}
\usepackage{mathptmx}  
\usepackage{float}
\usepackage[nodots]{numcompress}

\begin{document}

\markboth{A.~Mushtaq, A.~Kv{\ae}rn{\o}, and  K.~Olaussen}
{Higher order geometric integrators for a class  of Hamiltonian systems}

%
\catchline{}{}{}{}{}

\title{HIGHER ORDER GEOMETRIC INTEGRATORS FOR A CLASS OF
  HAMILTONIAN SYSTEMS}

\author{ASIF MUSHTAQ}
\address{Institutt for matematiske fag, NTNU, \\ N-7034, TRONDHEIM,
  NORWAY\,\\
\email{asifm@math.ntnu.no} }

\author{ANNE KV{\AE}RN{\O}}
\address{Institutt for matematiske fag, NTNU, \\ N-7034, TRONDHEIM,
  NORWAY\,\\
\email{anne.kvarno@math.ntnu.no} }

\author{K{\AA}RE OLAUSSEN}
\address{Institutt for fysikk, NTNU, \\ N-7491, TRONDHEIM,
  NORWAY\,\\
\email{kare.olaussen@ntnu.no} }

\maketitle

\begin{history}
\received{(Day Month Year)}
\revised{(Day Month Year)}
\end{history}

\begin{abstract}
We discuss systematic extensions of the standard (St{\"o}rmer-Verlet)
method for integrating the differential equations of Hamiltonian mechanics.
Our extensions preserve the symplectic geometry exactly, as well as all
N{\"o}ther conservation laws caused by joint symmetries of 
the kinetic and potential energies (like angular momentum in rotation 
invariant systems). These extentions increase the
accuracy of the integrator, which for the St{\"o}rmer-Verlet method is of order $\tau^2$ 
for a timestep of length $\tau$, to higher orders in $\tau$. The
schemes presented have, in contrast to most previous proposals,
all intermediate timesteps real and positive. They increase the relative accuracy
to order $\tau^{N}$ (for $N=4$, $6$, and $8$) for a large class of Hamiltonian systems. 

\end{abstract}


\section{Introduction}
\label{intro}
Numerical methods for integrating dynamical systems forward in time
will usually introduce errors, sometimes leading to results which
are even qualitatively wrong. However, there exists a class of methods,
geometric integrators, which aims at preservation of
the basic geometric properties of the system. A broad discussion of
such methods can be found in the book by Hairer \emph{et. al.}~\cite{HLWG}.
For many dynamical systems the
relevant geometry is the \emph{symplectic structure}
of phase space. For a phase space of dimension $2\cal{N}$ this structure can f.i.~be defined by the \emph{Poisson bracket}
\begin{equation}
     \left\{ \mathcal{A}(\bm{q}, \bm{p}), 
       \mathcal{B}(\bm{q}, \bm{p}) \right\} \equiv  \sum^{\cal{N}}_{a=1} 
     \left(\frac{\partial \mathcal{A}}{\partial q^a} \frac{\partial \mathcal{B}}{\partial p_a} - 
     \frac{\partial \mathcal{B}}{\partial q^a} \frac{\partial \mathcal{A}}{\partial p_a}\right)(\bm{q},\bm{p}),
\end{equation}
or the differential two-form,
\begin{equation}
         \Omega  = \sum^{\cal{N}}_{a=1} dq^a \bm{\wedge} dp_a.
\end{equation}
Another class of geometric constraints are those imposed by
conservation laws due to continuous (N{\"o}ther) symmetries.
An interesting field of applications for geometric integrators
are (classical limits of) \emph{minisuperspace} models of cosmology;
for a good introduction see the review by Capozziello 
\emph{et. al.}~\cite{CapozzielloDeLaurentisOdintsov}. Some of us
are in the process of analysing questions from this field by use
of the algorithms described in this paper.

The Hamilton equations of motion constitute a system of
ordinary first order differential equations,
\begin{align}
   \dot{q}^a = \frac{\partial H}{\partial p_a},\quad
   \dot{p}_a = -\frac{\partial H}{\partial q^a}, \quad a=1,\ldots,\cal{N},
   \label{HamiltonEquation}
\end{align}
where the overdot $\dot{\ }$ denotes differentiation with respect to time $t$, and
$H = H(\bm{q}, \bm{p})$.  They can be viewed as the characteristic
equations of the partial differential equation
\begin{equation}
    \frac{\partial}{\partial t} \rho(\bm{q}, \bm{p}; t) =
    {\cal L}_{H}\,\rho(\bm{q}, \bm{p}; t),
    \label{HamiltonianFlow}
\end{equation}
with ${\cal L}_H$ the first order differential operator,
\begin{equation}
    {\cal L}_{H}\;\bm{\cdot}\; = \sum_{a=1}^{\cal{N}}
    \left(\frac{\partial H}{\partial p_a}\frac{\partial}{\partial q^a}
    -\frac{\partial H}{\partial q^a}\frac{\partial}{\partial p_a}\right) \;\bm{\cdot}\; \equiv \left\{ \;\bm{\cdot}\;, H \right\},
\end{equation}
generating a flow on phase space.
If $H$ does not depend ex\-pli\-citly on $t$, a formal solution
of (\ref{HamiltonianFlow}) is
\begin{equation}
    \rho(\bm{q}, \bm{p}; t) = 
    \text{e}^{t{\cal L}_H} \rho(\bm{q}, \bm{p}; 0).
\end{equation}
In most cases this expression remains just formal, but one may
often split the Hamiltonian into two parts,
$H = H_A + H_B$, corresponding to a splitting
${\cal L}_H = {\cal L}_A + {\cal L}_B$, such that the flows generated by
${\cal L}_A$ and ${\cal L}_B$ are separately integrable. One may
then use the Cambell-Baker-Hausdorff formula to approximate
the flow generated by ${\cal L}_H$. By doing this in a symmetric way one
obtains the Strang~\cite{Strang} splitting formula,
\begin{align}
     &\exp\left({\textstyle {\frac{1}{2}}\tau{\cal L}_B}\right)\,
     \exp\left({\tau{\cal
           L}_A}\right)\,\exp\left({\textstyle {\frac{1}{2}}\tau{\cal L}_B}\right)
    = \nonumber\\ &  
    \exp\left({\tau(\mathcal{L}_A + \mathcal{L}_B ) + 
       {\textstyle \frac{1}{24}}\tau^3 \left[2{\cal L}_A
         + {\cal L}_B,\left[{\cal L}_A,{\cal L}_B\right]\right]+\cdots}\right),
     \label{CambellBakerHousdorff}
\end{align}
which shows that time stepping this expression with a timestep $\tau$
provides an approximation with relative accuracy of order $\tau^2$,
exactly preserving the symplectic property of the flow. 

The left hand side of equation~\eqref{CambellBakerHousdorff} corresponds
to the symplectic splitting scheme of solving
\begin{subequations}
\begin{align}
       \dot{q}^a &= \frac{\partial H_B}{\partial p_a} ,\quad
       \dot{p}_a  = -\frac{\partial H_B}{\partial q^a},\quad\text{for
         a timestep $\frac{1}{2}\tau$,}\label{FirstKick}\\
       \dot{q}^a &= \frac{\partial H_A}{\partial p_a} ,\quad
       \dot{p}_a  = -\frac{\partial H_A}{\partial q^a},\quad\text{for
         a timestep $\tau$,}\label{Move}\\
       \dot{q}^a &= \frac{\partial H_B}{\partial p_a} ,\quad
       \dot{p}_a  = -\frac{\partial H_B}{\partial q^a},\quad\text{for
         a timestep $\frac{1}{2}\tau$,\label{SecondKick}}
\end{align}
\label{SymplecticSplitting}
\end{subequations}

\noindent
and iterating.
Here the last part of one iteration may be combined with the
first part of the next, unless one wants to register the state of
the system at intermediate times.

A rather common situation
is that $H(\bm{q}, \bm{p}) = T(\bm{p}) + V(\bm{q})$, where one may
choose the splitting $H_A = T$ and 
$H_B = V$ (or $H_A = V$ and $H_B = T$). 
In quantum mechanics such a splitting scheme is usually associated with
the Trotter product formula,
\begin{equation}
    \text{e}^{-\text{i}(T+V) t} =
    \mathop{\lim}_{n\to\infty}\left(\text{e}^{-\text{i}Tt/n}
    \text{e}^{-\text{i} Vt/n}\right)^n,
\end{equation}  
but there more often used to motivate the path
integral formulation than as a numerical approximation method. 
For a more common language we refer to steps \eqref{FirstKick} and
\eqref{SecondKick} as \emph{kicks}, since they (when $H_B = V$)
change momentum but not position, and step \eqref{Move} as a \emph{move},
since it change position but not momentum.

The integration scheme defined by the equations \eqref{SymplecticSplitting},
although symplectic, does not quite reproduce the flow generated by $H$, but instead one
generated by
\begin{equation}
     H'(\bm{q},\bm{p}) = H(\bm{q}, \bm{p}) - \tau^2 H_2(\bm{q},\bm{p})
     + \mathcal{O}(\tau^4),
     \label{StormerVerletHamiltonian}
\end{equation}
for some Hamiltonian $H_2$ (computed in section \ref{nonLinearSystems}).
For this reason the motion will not be exactly confined to
the constant energy surface $\mathcal{S}_H(E)$ of $H$, but to the constant energy surface
$\mathcal{S}_{H'}(E')$ of $H'$,
which for small enough $\tau^2$ (and regular Hamiltonians) will
lie close to $\mathcal{S}_H(E)$. Hence, the energy error of the generated solution will
stay bounded regardless of how long we integrate the equations, being geometrically
constrained by the maximal distance between $\mathcal{S}_H(E)$ and $\mathcal{S}_{H'}(E')$.
In addition to symplecticity, this constraint is another reason why the scheme defined
by \eqref{SymplecticSplitting} is qualitatively robust.

Nevertheless, the generated and exact solutions may still drift apart with time, 
since two points on (essentially) the same solution surface may still be quite
far apart. Hence, there are good reasons to search for more accurate schemes
which maintain the attractive geometric properties of equation~\eqref{SymplecticSplitting}.
Our approach is to stick to the three-step algorithm of~\eqref{SymplecticSplitting}, but
with modified expressions for the Hamiltonians $H_A$ and $H_B$. For instance,
finding the error term of equation~\eqref{StormerVerletHamiltonian} to
be $H_2 = T_2 + V_2$,
we can eliminate the order-$\tau^2$ error by using an effective kinetic energy 
$T_{\text{eff}} = T + \tau^2 T_2 + \mathcal{O}(\tau^4)$, 
and an effective potential energy
$V_{\text{eff}} = V + \tau^2 V_2 + \mathcal{O}(\tau^4)$ in the splitting scheme.

There have been several approaches to construct integration schemes which
are of higher order in $\tau$, while maintaining exact symplecticity
of the evolution.
Accessible reviews of such approaches have f.i. been given by Yoshida~\cite{YoshidaReview},
McLachan \emph{et.~al.}~\cite{RIM},
and Blanes \emph{et.al.}~\cite{BlanesCasasMurua}.
Neri~\cite{Neri}  has provided the general idea to construct
symplectic integrators for Hamiltonian systems.
Forest and Ruth ~\cite{FR} discussed an explicit fourth
order method for the integration of Hamiltonian equations for the
simplest non-trivial case.
Suzuki~\cite{Suzuki} presented the idea of how recursive construction
of successive approximants may be extended to other methods. 

Many of the higher order symplectic splitting methods involve an extension of
equation~(\ref{CambellBakerHousdorff}) to an expression of the form
\begin{equation}
     \prod_{i=1}^k\, \exp\left({c_i\tau{\cal L}_B}\right)\,\exp\left({d_i\tau{\cal L}_A}\right) = 
     \exp\left({\tau ({\cal L}_A +{\cal L}_B)} + {\cal O}(\tau^{\ell})\right),
     \label{GeneralIteratedComposition}
\end{equation}
as discussed by Yoshida~\cite{Yoshida}. It was noted that if one uses a symmetric integrator, such that
\begin{equation*}
    S(\tau) = \exp\left({\tau({\cal L}_A +{\cal L}_B)} + \tau^{\ell}\,R_{\ell} + {\cal O}(\tau^{\ell+2})\right)
\end{equation*}
for some generator $R_{\ell}$, then
\begin{equation*}
    S(x_1\tau)\, S(x_0\tau)\,S(x_1\tau) = 
    \exp\left({(x_0 + 2x_1)\tau({\cal L}_A +{\cal L}_B)} + 
   (x^{\ell}_0 + 2 x^{\ell}_1)\tau^{\ell}\,R_{\ell}  +{\cal O}(\tau^{{\ell}+2})\right).
\end{equation*}
Hence, by choosing
\begin{subequations}
\label{ImprovementConditions}
\begin{align}
  x_0 + 2x_1 &= 1,    \label{AddToOne}\\
  x^{\ell}_0 + 2x^{\ell}_1 &= 0, \label{AddToZero}
\end{align}
\end{subequations}
one increases the order of the scheme by two or more. However,
equations~(\ref{ImprovementConditions})
have real solutions only if either $x_0$ or $x_1$ is negative.
In fact, it has been
proven (cf.~Sheng \cite{Sheng}, Suzuki \cite{SuzukiJMP}, Goldman and Kaper \cite{GoldmanKaper}, 
Blanes and Casas \cite{BlanesCasas}) that
all schemes of the form~(\ref{GeneralIteratedComposition})
require at least one $c_i<0$, and at least one $d_i<0$. For equations invariant under
time reversal, which is often the case for Hamiltonian systems,
this may not be a big obstacle (although it may seem like an inefficient
way of integrating equations forward in time). 

Worse, if one wants to use the same code to solve
parabolic equations\footnote{Like solving a heat type equation,
$-\partial_t u(t,\bm{x}) = \left[-{\Delta} + V(\bm{x})\right] u(t,\bm{x})$,
instead of a Schr{\"o}dinger type equation,
$\text{i}\partial_t \psi(t,\bm{x}) = \left[-{\Delta} + V(\bm{x})\right] \psi(t,\bm{x})$,
obtained by formally replacing $t$ with $-\text{i}\,t$ in the Schr{\"o}dinger equation.}
negative timesteps may have a disastrous effect on numerical stability due to exponentially
growing errors.
Castella \emph{et. al.}~\cite{CastillaChartierDescombesVilmart} have proposed 
to use \emph{complex} solutions of equations~(\ref{AddToOne}, \ref{AddToZero}).
It is possible to find solutions where all timesteps have a positive real part.
This can stabilize the scheme, but at the cost of working with complex quantities.

In this paper we investigate a different approach,
based on our~\cite{AsifAnneKare} observation that the operators 
${\cal \tilde{L}}_A$, ${\cal \tilde{L}}_B$ of each step of a
splitting scheme don't need to be \emph{exactly} the same as those
in the sum ${\cal L}_H = {\cal L}_A + {\cal L}_B$. Instead, our approach
is to construct $\tau$-dependent operators 
${\cal \tilde{L}}_A$, ${\cal \tilde{L}}_B$ such that
\begin{equation}
    \exp\left({{\textstyle \frac{1}{2}}\tau{\cal{\tilde L}}_B}\right)\,
    \exp\left({\tau {\cal \tilde{L}}_A}\right)\,
    \exp\left({{\textstyle \frac{1}{2}}\tau \tilde{\cal L}_B}\right) = 
    \exp\left({\tau({\cal L}_A + {\cal L}_B) + {\cal O}(\tau^k)}\right).
\end{equation}
For Hamiltonians of the form, 
\begin{equation}
   H(\bm{q},\bm{p}) = \frac{1}{2}\,\bm{p}^T M \bm{p} + V(\bm{q}),
   \label{SeparableHamiltonian}
\end{equation}
with $M$ a symmetric positive definite matrix (the inverse mass matrix),
we have constructed an explicit series expansion
\begin{equation}
      {\cal\tilde{L}}_A = {\cal L}_A + \sum_{k=1}^{(N-2)/2}
      \tau^{2k}\,{\cal L}^{(k)}_A,\quad
      {\cal\tilde{L}}_B = {\cal L}_B + \sum_{k=1}^{(N-2)/2} \tau^{2k}\,{\cal L}^{(k)}_B.
\end{equation}
up to $N=8$. This can be used in schemes with global error of order
$\tau^N$, for $N\in \left\{2, 4, 6, 8\right\}$.
We denote the order of these schemes by $N$.
Since the operators $\tilde{\cal L}_X$ (with $X = A\text{ or }B$)
generate flows $\exp\left(\tau \tilde{\cal L}_X\right)$ which
are modifications of those generated by ${\cal L}_X$,
we refer to such flows as \emph{modified integrators}.
Chartier \emph{et.~al.}~\cite{ChartierHairerVilmart}
have labeled such schemes as \emph{modified differential
equations}. 
The $N$'th order scheme is constructed to generate the same flow as 
\begin{equation}
   {\cal L}' = {\cal L} - 
   \tau^{N} \left(\tilde{\cal L}^{(N/2)}_A + \tilde{\cal L}^{(N/2)}_B\right) + {\cal O}(\tau^{N+2}),
\end{equation}
after each complete timestep.
I.e., we use modified integrators to generate the unmodified flow
better. Note that knowledge of the Hamiltonians corresponding to the
operators $\tilde{\mathcal{L}}^{(N/2)}_X$ also provide an explicit estimate of
the leading error of the scheme. 

One possible restriction on the class of available splitting schemes
is the requirement that both of the flows
$\exp\left(\frac{1}{2}\tau\tilde{\cal L}_B\right)$ 
and $\exp\left(\tau\tilde{\cal L}_A\right)$ must
be \emph{ex\-plic\-itly integrable}.  We have relaxed this requirement
by demanding both flows to be \emph{efficiently computable}:
I.e., it must be possible to integrate each short timestep numerically
sufficiently fast, while preserving the symplectic structure to
sufficient numerical precision.
This question arises for our method because $T_{\text{eff}}$ in general will
depend on both $\bm{p}$ and $\bm{q}$. This means that 
the step defined by equation \eqref{Move} will in general not be explicitly integrable.
Instead, to integrate such steps we construct a generating function $G$
for a transformation
\begin{equation}
  {\cal G}: \left\{\bm{q}(t),\bm{p}(t)\right\}
  \to \left\{ \bm{q}(t+\tau),\bm{p}(t+\tau) \right\},
\end{equation}
which solves equation~\eqref{Move} to the required accuracy in $\tau$,
and is exactly symplectic. The change in momentum $\bm{p}$ (of
order $\tau^3$ --- i.e.~a gentle \emph{push}) is then defined
through an implicit equation, while the change in
position $\bm{q}$ continues to be explicit.


The rest of this paper is organized as follows: In section~\ref{linearSystems}
we demonstrate the basic idea of the proposed methods on linear systems. Next we
develop the general theory, valid for separable Hamiltonians~(\ref{SeparableHamiltonian}),
in section~\ref{nonLinearSystems}. Here we construct the operators $\tilde{\cal L}_X$ (for $X=A, B$) 
explicitly. Or more precisely, we calculate the contributions $T_{2k}$ to the corresponding $T_{\text{eff}}$,
and the contributions $V_{2k}$ to the corresponding $V_{\text{eff}}$, cf.~equation \eqref{ModifiedHamiltonians}.
We focus our discussion on the numerical implementations
in section~\ref{numericalResults}, together with test-case investigations
of how the methods work in practice.
We have tested these methods on anharmonic oscillators and
Fermi-Pasta-Ulam-Tsingou type problems (named as suggested 
by Dauxois~\cite{Dauxois}). We conclude the paper with some brief
remarks in section~\ref{concludingRemarks}.

\section{Linear Systems} \label{linearSystems}

\subsection{Single Harmonic Oscillator}

For a simple illustration of our idea consider the Hamiltonian
\begin{align}
   H({q},{p})=\frac{1}{2}\left( p^2 + q^2\right),
\end{align}
 whose exact evolution over a time interval $\tau$ is
\begin{equation}
     \begin{pmatrix} q_{\text{e}} \\[0.4ex]  p^{\text{e}} \end{pmatrix} =
    \begin{pmatrix*}[r] \cos \tau &  \sin \tau \\  - \sin \tau &
  \cos\tau \end{pmatrix*}   
   \begin{pmatrix} q\\p\end{pmatrix}.
\label{H_0a}
\end{equation}
Compare this with the process of first evolving the system with the
(\emph{kick}) Hamilton\-ian 
\(
  H_{B}= {\frac{1}{2} k q^2}
\)
for a time $\frac{1}{2}\tau$,
followed by an evolution with the (\emph{move}) Hamilton\-ian
\(
  H_{A}= {\frac{1}{2} m p^2}
\)
for a time $\tau$, and ending with an evolution with
$H_{B}$ for a time $\frac{1}{2}\tau$.
One such combination (one complete timestep) gives
\begin{equation}
    \begin{pmatrix}q_{\text{s}} \\ p^{\text{s}}\end{pmatrix} =
    \begin{pmatrix*}[c] 1-\frac{1}{2}m k\tau^2 & m \tau \\[0.4ex]
      -(1-\frac{1}{4} km \tau^ 2) k\tau &  1-\frac{1}{2} k m\tau^2 
    \end{pmatrix*}
    \begin{pmatrix}q \\ p \end{pmatrix}.
\end{equation}
We note that by choosing
\begin{subequations}
\label{HarmonicOscillatorCorrection}
\begin{align}
    m &= \frac{\sin\tau}{\tau} = 1 - \frac{1}{6}\tau^2 +
    \frac{1}{120}\tau^ 4 - \frac{1}{5040} \tau^6 + \cdots,\\
    k &= \frac{2}{\tau} \tan \frac{\tau}{2} =
    1+\frac{1}{12}\tau^2+\frac{1}{120}\tau^4 + \frac{17}{20160}\tau^6
    + \cdots,
\end{align}
\end{subequations}
the exact evolution is reproduced, provided $0 < \tau < \pi$.
If we interchange the r{\^o}les of
$H_{A}$ and $H_{B}$ one combination instead gives
\begin{equation}
    \begin{pmatrix}q_{\text{s}} \\[0.4ex] p^{\text{s}}\end{pmatrix} =
    \begin{pmatrix*}[c] 1-\frac{1}{2}{m} {k}\tau^2 & 
      (1-\frac{1}{4}{m}{k}){m} \tau \\[0.5ex]
      -{k}\tau &  1-\frac{1}{2} {k}{m}\tau^2 
    \end{pmatrix*}
    \begin{pmatrix}q \\ p \end{pmatrix},
\end{equation}
which becomes exact if we choose
\begin{equation}
  {m} = \frac{2}{\tau}\tan\frac{\tau}{2},\quad
  {k} = \frac{\sin\tau}{\tau},
  \label{ModifiedParameters}
\end{equation}
again provided the timestep is restricted to the interval $0 < \tau < \pi$.

\subsection{Higher-dimensional linear systems}

It should be clear that this idea works for systems of harmonic
oscillators in general, i.e. for quadratic Hamiltonians of the form
\begin{equation}
   H(\bm{q},\bm{p}) = \frac{1}{2} \left( \bm{p}^T M \bm{p} + \bm{q}^T
     K \bm{q} \right),
   \label{GeneralHamiltonian}
\end{equation}
where $M$ and $K$ are symmetric matrices.
For a choosen splitting scheme and step 
interval $\tau$ there are always modified matrices,
\begin{subequations}
\begin{align}
     M_\tau &= M - \frac{\tau^2}{6}  M\,(K\,M) + \frac{\tau^4}{120} M\,(K\,M)^2 
     - \frac{\tau^6}{5040} M\,(K\,M)^3 + {\cal O}(\tau^8),\\
     K_\tau &= K + \frac{\tau^2}{12} (K\,M)\,K + \frac{\tau^4}{120} (K\,M)^2\,K 
     + \frac{17\tau^6}{20160} (K\,M)^3\,K + {\cal O}(\tau^8),
\end{align}
\end{subequations}
generating a \emph{kick-move-kick} flow which reproduces the exact
one up to corrections of order $\tau^8$.
It should be obvious how this can be extended to arbitrary order
in $\tau^2$, with coefficients taken from the expansions in
equation (\ref{HarmonicOscillatorCorrection}).
In principle this can be used to reproduce the exact flow,
provided $\tau$ is not too large. I.e., $0 < \tau < \pi/\omega_{\text{max}}$,
where $\omega_{\text{max}}$ is the largest angular frequency of the system.

\section{General potentials} \label{nonLinearSystems}

For a more general treatment we consider Hamiltonians of the form
\begin{equation}
    H(\bm{q}, \bm{p}) = \frac{1}{2} \bm{p}^T M \bm{p}  + V(\bm{q}).
\end{equation}
A series solution of the Hamilton equations in powers of $\tau$ is
\begin{subequations}
\begin{align}
     q^a_{\text{e}} &= q^a + p^a\tau - \frac{1}{2}\partial^a V \tau^2 
     -\frac{1}{6}\partial^a (DV) \tau^3+{\cal O}(\tau^4),\\
     p^{\text{e}}_a &= p_a - \partial_a V \tau  -
     \frac{1}{2}\partial_a (DV) \tau^2 + \partial_a 
\left( \frac{1}{12} \bar{D}V -\frac{1}{6} D^2V\right) \tau^3 + {\cal O}(\tau^4).
\end{align}
\end{subequations}
Here we have introduced notation to compactify expressions,
\begin{align}
     &\partial_a \equiv \frac{\partial}{\partial q^a},\quad 
     \partial^a \equiv M^{ab}\partial_b,\quad
     p^a \equiv M^{ab} p_b,\quad D \equiv p_a \partial^a,\quad
     \bar{D} \equiv (\partial_a V)\partial^a,\nonumber
\end{align}
where we employ the {\em Einstein summation convention\/}: An
index which occur twice, once in lower position and once in upper
position, are implicitly summed over all its available values. I.e.,
\begin{equation} 
  M^{ab} \partial_b \equiv \sum_b M^{ab} \partial_b.
\end{equation}
We will generally use the matrix $M$ to rise an index from lower to upper position.

If we instead use a splitting method to generate the flow,
with generators $H_B = V(\bm{q})$ and $H_A = \frac{1}{2} \bm{p}^T M \bm{p} \equiv T$
(i.e., a {\em kick-move-kick\/} scheme), we obtain
\begin{subequations}
\begin{align}
    q_{\text{s}}^a &= q^a + p^a\tau - \frac{1}{2}\partial^a V \tau^2  +{\cal O}(\tau^4),\\
    p^{\text{s}}_a &=p_a - \partial_a V\tau - \frac{1}{2}\partial_a (DV) \tau^2 
     +\partial_a \left( \frac{1}{8}\bar{D}V  - \frac{1}{4}D^2V\right) \tau^3 + {\cal O}(\tau^4).
\end{align}
\end{subequations}
As expected the result differs from the exact result in the third order.
However,  the difference can be corrected by modifiying the generators,
$H_A \rightarrow T + \tau^2\,T_2$ and $H_B \rightarrow V + \tau^2\,V_2$, with
\begin{equation}
    T_2 = -\frac{1}{12} D^2V,\quad
    V_2 = \frac{1}{24} \bar{D}V.
    \label{SecondOrderCorrections}
\end{equation}
Specialized to a one-dimensional system with potential $V = \frac{1}{2}q^2$
this agrees with equation~(\ref{HarmonicOscillatorCorrection}). With
this correction the {\em kick-move-kick\/} splitting scheme agrees
with the exact solution to fourth order in $\tau$, but differ
in the $\tau^5$-terms. We may again correct the difference by introducing
fourth order generators, $H_A \rightarrow T + \tau^2\,T_2 + \tau^4\,T_4$ and
$H_B \rightarrow V + \tau^2\,V_2 + \tau^4\,V_4$, with
\begin{align}
    T_4 =\frac{1}{720}\left( D^4  - 9 \bar{D} D^2 + 3 D\bar{D} D\right) V,
    \quad V_4=\frac{1}{480} \bar{D}^2 V.
    \label{FourthOrderCorrections}
\end{align}
Specialized to a one-dimensional system with potential $V = \frac{1}{2}q^2$
this agrees with equation~(\ref{HarmonicOscillatorCorrection}). 
With this correction the {\em kick-move-kick\/} splitting scheme agrees
with the exact solution to sixth order in $\tau$, but differ
in the $\tau^7$-terms. We finally correct this difference by introducing
sixth order generators, $H_A \rightarrow T + \tau^2\,T_2 + \tau^4\,T_4 + \tau^6\,T_6$
and $H_B \rightarrow V + \tau^2\,V_2 + \tau^4\,V_4 + \tau^6\,V_6$,
with
\begin{subequations}
\label{SixthOrderCorrections}
\begin{align}
    T_6=& -\frac{1}{60480}\left(2\, D^6-40\, \bar{D}D^4 +46\,
      D\bar{D} D^3 -15\,D^2 \bar{D}D^2 +54\,\bar{D}^2 D^2\right. \nonumber\\  
       & - 9\,\bar{D}D\bar{D}D - 42\,D\bar{D}^2 D\left. +12\, D^2\bar{D}^2\right)V\\
    V_6 &= \frac{1}{161280}\left( 17\, \bar{D}^3 - 10\,\bar{D}_3\right)V,
\end{align}
\end{subequations}
where we have introduced 
\begin{equation}
 \bar{D}_3 \equiv (\partial_aV)(\partial_b V)(\partial_c
 V) \partial^a \partial^b \partial^c.
\end{equation}
Specialized to a one-dimensional system with potential $V = \frac{1}{2}q^2$
this agrees with equation~(\ref{HarmonicOscillatorCorrection}). 
With this correction the {\em kick-move-kick\/} splitting scheme agrees
with the exact solution to eight order in $\tau$, but differ
in the $\tau^9$-terms. The process may be continued to higher orders
in $\tau$, 
\begin{align}
  H_A \rightarrow T + \sum_{k\ge1} \tau^{2k}\,T_{2k},\qquad H_B \rightarrow 
  V + \sum_{k\ge1} \tau^{2k}\,V_{2k}.\label{ModifiedHamiltonians}
\end{align}

To keep track of the algebraic expressions which occured during 
the calculations above, we have represented them graphically
in terms of bi-colored tree-diagrams. I.e., these calculations are related
to ``rooted-tree-type'' theories. Our tree-diagrams describing $T_{2k}$ and $V_{2k}$,
and the generating functions $G_k$ below, are unrooted (the derivatives of these
scalar functions can be represented by rooted trees).
It is fairly straightforward to find
the general structure of the order $\tau^{N}$ correction terms,
but more laborious to compute
the rational coefficients multiplying each term.
They are simplest found by considering enough
special cases for a unique determination. After the
explicit expressions (\ref{SecondOrderCorrections}, \ref{FourthOrderCorrections})
were found we verified them manually for a general
Hamiltonian (\ref{GeneralHamiltonian}) using graphical calculations.
The explicit expressions (\ref{SixthOrderCorrections}) has been checked
against a general Hamiltonian (\ref{GeneralHamiltonian}) acting on a
four-dimensional phase space (i.e., with two-dimensional $\bm{q}$ and
$\bm{p}$) by use of a computer algebra program.

\section{Solving the {\em move\/} steps}

Addition of extra potential terms 
$V \rightarrow V_{\text{eff}} \equiv V + \tau^2\,V_2 + \tau^4\,V_4 + \dots$, 
is in principle unproblematic for solution of the {\em kick\/} steps. The
equations,
\begin{equation}
  \dot{q}^a = 0,\quad
  \dot{p}_a = -\partial_a V_{\text{eff}}(\bm{q}),
\end{equation}
can still be integrated exactly, preserving the symplectic structure.
The situation is different for the kinectic term
$T \rightarrow T_{\text{eff}} \equiv T + \tau^2\,T_2 + \tau^4\,T_4 + \cdots$,
since it now leads to equations
\begin{equation}
  \dot{q}^a = \frac{\partial }{\partial p_a}T_{\text{eff}}(\bm{q},\bm{p}),\quad
  \dot{p}_a = -\partial_a T_{\text{eff}}(\bm{q},\bm{p}),
  \label{MoveSteps}
\end{equation}
which are no longer straightforward to integrate exactly.
Although the problematic terms are small
one should make sure that the {\em move\/} steps preserve
the symplectic structure. Let $\bm{q}, \bm{p}$
denote the positions and momenta just before the {\em move\/} step, and $\bm{Q}, \bm{P}$
the positions and momenta just after. The relation between 
$\bm{q}, \bm{p}$ and $\bm{Q}, \bm{P}$ can be expressed in terms of a 
generating  function (cf. Goldstein~\cite{Goldstein}, Arnold~\cite{Arnold}),
\begin{equation}
   G(\bm{q}, \bm{P};\tau) = q^a P_a + \Delta G(\bm{q}, \bm{P};\tau),
\end{equation}
such that the transformation
\begin{equation}
   Q^a = \frac{\partial G}{\partial P_a},\qquad
   p_a= \frac{\partial G}{\partial q^a} = P_a + 
   \frac{\partial \Delta G}{\partial q^a}.\label{CanonicalTransformation}
\end{equation}
preserves the symplectic structure exactly. However,
note that the relation between $\bm{p}$ and $\bm{P}$ in general
is a nonlinear equation of the form
\begin{equation}
   P_a = p_a - \frac{\partial}{\partial q^a}\Delta G(\bm{q},
   \bm{P};\tau),
   \label{PushStep}
\end{equation}
where the second term on the right is of order $\tau^3$ or higher. We solve
this equation by iteration. With $\bm{P}^{(0)} = \bm{p}$,
\begin{equation}
    P^{(n+1)}_a = p_a - \frac{\partial}{\partial q^a}\Delta G(\bm{q}, \bm{P}^{(n)};\tau).
\end{equation}
Writing $\bm{P}^{(n)} = \bm{P} +\bm{\Delta P}^{(n)}$, with $\bm{P}$ the exact solution,
we find to first order in $\bm{\Delta P}$ that
\begin{equation}
   {\Delta P}^{(n+1)}_a =  -\frac{\partial^2}{\partial q^a\partial
     P_b}\Delta G(\bm{q}, 
   \bm{P};\tau)\,{\Delta P}^{(n)}_b
   \equiv - \Delta G_a^{\phantom{a}b}\,{\Delta P}^{(n)}_b.
\end{equation}
Let $\lambda \sim \tau^3$ be the eigenvalue of $\Delta G_a^{\phantom{a}b}$ with largest magnitude.
Then the iteration converges exponentially fast towards the exact solution,
with $\bm{\Delta P}^{(n)}$ decaying like $\lambda^n \sim \tau^{3n}$. Since it is most to gain
by a higher order method when the timestep $\tau$ is small, we
assume $\lambda$ to be small in all cases of practical relevance.
Our experience is that the iteration scheme is robust,
with 3--4 iterations been sufficient for computations
to double precision accuracy. It is important that \eqref{PushStep} is
solved to sufficient accuracy; otherwise the evolution fails to be
(sufficiently) symplectic.


Some of our theoretical results have already been given in the literature. 
The generating function formalism has been used earlier by Feng~\cite{Feng}
and Feng \emph{et.al} ~\cite{FengWuQinWang}
to construct canonical difference schemes  (see also Channell
and Scovel~\cite{ChannelScovel}, Stuchi~\cite{ Stuchi}). They give
the result~(\ref{SecondOrderCorrections}), but the actual solution of
the resulting implicit equations are not discussed.
One can construct a generating function for the full symplectic
evolution over a timestep $\tau$, without combination with a splitting method.
However, in that case the resulting nonlinear equations would
be more time consuming and/or difficult to solve by direct iteration.

We now explicitly construct
$G$ so that the {\em move\/} step is reproduced to sufficient accuracy.
Consider first the case when $H_A = T$. The choice
$G = q^a P_a + \frac{1}{2} P^a P_a \tau$ gives
\begin{equation}
     Q^a = q^a + P_a\,\tau,\quad p_a = P_a,
     \label{SimpleKick}
\end{equation}
which is the correct relation. Now add the $\tau^2\,T_2$-term to the {\em  move\/} step.
The exact solution of equation (\ref{MoveSteps}) becomes 
\begin{subequations}
\label{ExactMove}
\begin{align}
    Q^a &= q^a + p_a\,\tau - {\textstyle \frac{1}{6}} \partial^a DV\,\tau^3 
    -{\textstyle \frac{1}{24}}\partial^a D^2V\, \tau^4 + {\cal O}(\tau^5),\\
    P_a &= p_a + {\textstyle \frac{1}{12}} \partial_a \, D^2 V\,\tau^3 
    + {\textstyle \frac{1}{24}}\partial_a D^3 V\,\tau^4 + {\cal O}(\tau^5).
\end{align}
\end{subequations}
Compare this with the result of changing
\begin{equation}
  G \rightarrow G -{\textstyle \frac{1}{12}} {\cal D}^2\, V \tau^3 
  - {\textstyle \frac{1}{24}} {\cal D}^3V\, \tau^4,
\end{equation}
where ${\cal D} \equiv P_a \partial^a$. The solution of equation~(\ref{CanonicalTransformation})
change from the relations~(\ref{SimpleKick}) to
\begin{subequations}
\begin{align}
  Q^a &= q^a + P^a\tau - {\textstyle \frac{1}{6}}\partial^a {\cal D}V \tau^3 
  -{\textstyle \frac{1}{8}}\partial^a {\cal D}^2V\,\tau^4 + {\cal O}(\tau^5),\label{Qequation}\\
  p_a &= P_a - {\textstyle \frac{1}{12}} \partial_a {\cal D}^2 V \tau^3 
  - {\textstyle \frac{1}{24}}\partial_a{\cal D}^3 V \tau^4 + {\cal O}(\tau^5). \label{Pequation}
\end{align}
\end{subequations}
Since ${\cal D}$ is  linear in $\bm{P}$, 
equation~(\ref{Pequation})
constitute a system of third order algebraic equations which in general
must be solved numerically. This should usually be a fast process for
small $\tau$. An exact solution of this equation is required to preserve
the symplectic structure, but  this solution should also agree with
the exact solution of (\ref{MoveSteps}) to order $\tau^4$. This may be
verified by perturbation expansion in $\tau$. A perturbative solution
of equation~(\ref{Pequation}) is 
\begin{equation*}
    P_a = p_a + {\textstyle \frac{1}{12}} \partial_a D^2 V \tau^3 
    + {\textstyle \frac{1}{24}} \partial_a D^3 V \tau^4 + {\cal O}(\tau^5),
\end{equation*}
which inserted into (\ref{Qequation}) reproduces
the full solution~(\ref{ExactMove}) to order $\tau^4$.

This process can be systematically continued to higher orders.
We write the transformation function as
\begin{equation}
G(\tau) = \sum_{k=0}^\infty \tau^k\,G_k,
\end{equation}
and find the first terms in the expansion to be
{\footnotesize
\begin{align*}
    G_0& = q^a P_a,\quad G_1 = {\textstyle \frac{1}{2}} P^a P_a,\quad  G_2 =
    0,\quad  G_3 = -{\textstyle \frac{1}{12}}{\cal D}^2 V,\quad G_4 = 
    -{\textstyle \frac{1}{24}}{\cal D}^3 V,\nonumber\\
    G_5 &= -{\textstyle \frac{1}{240}}\left( 3\,{\cal D}^4+ 3\, \bar{D}{\cal D}^2
      -{\cal D} \bar{D} {\cal D} \right) V,\\
    G_6 &= -{\textstyle \frac{1}{720}} \left({2\,\cal D}^5  + 8\,\bar{D}{\cal D}^3
      - 5\, {\cal D}\bar{D}{\cal D}^2 \right) V,\nonumber\\
    G_7 &=-{\textstyle \frac{1}{20160}}\left({10\,\cal D}^6 + 10\,\bar{D}{\cal D}^4 + 90\,{\cal
        D}\bar{D} {\cal D}^3 - 75\,{\cal D}^2 \bar{D} {\cal D}^2\right.
    \nonumber
    \\ &\phantom{=-{\textstyle \frac{1}{20160}}}\left.\, +\,18\,\bar{D}^2{\cal D}^2 -3\, \bar{D}{\cal
         D}\bar{D}{\cal D} -14\,{\cal D}\bar{D}^2{\cal D} + 4\,
       {\cal D}^2\bar{D}^2 \right) V,\nonumber\\
    G_8 &= - {\textstyle \frac{1}{40320}} \left(3\,{\cal D}^7-87\,\bar{D}{\cal D}^5
      +231\,{\cal D}\bar{D}{\cal D}^4 -133\,{\cal D}^2\bar {D}{\cal
        D}^3+ 63\, \bar{D}^2 {\cal D}^3  -3\,{\cal D}\bar{D}^2{\cal D}^2
     -21\,{\cal D}^2 \bar{D}^2 {\cal D}   \right.\\ \nonumber
  &\phantom{= - {\textstyle \frac{1}{40320}}} \left.\, +\,4\,{\cal D}^3 \bar{D}^2-63
    \,\bar{D}{\cal D}\bar{D}{\cal D}^2  +25\,{\cal D}\bar{D}{\cal D}\bar{D}{\cal D}\right) V.\nonumber
\end{align*}}
Also in these calculations we represent the algebraic expressions by bi-colored tree diagrams,
to better visualise and understand their structure.
The possible graphical structures for $G_n$ is fairly simple to write down.
But it is quite laborious to find the rational coefficients
multiplying each graph. They are simplest found by considering enough
special cases for a unique determination. After that we have verified
the expressions up to $G_6$ manually using graphical
calculations, and $G_7$, $G_8$
against a general Hamiltonian (\ref{GeneralHamiltonian}) acting on a
four-dimensional phase space (i.e., with two-dimensional $\bm{q}$ and $\bm{p}$)
by use of a computer algebra program.

\section{Numerical results on nonlinear systems} \label{numericalResults}

\subsection{One-dimensional anharmonic oscillator}

\begin{figure}
\begin{center}
\includegraphics[clip, trim = 8ex 5.5ex 9ex 5ex, width=0.96\textwidth]{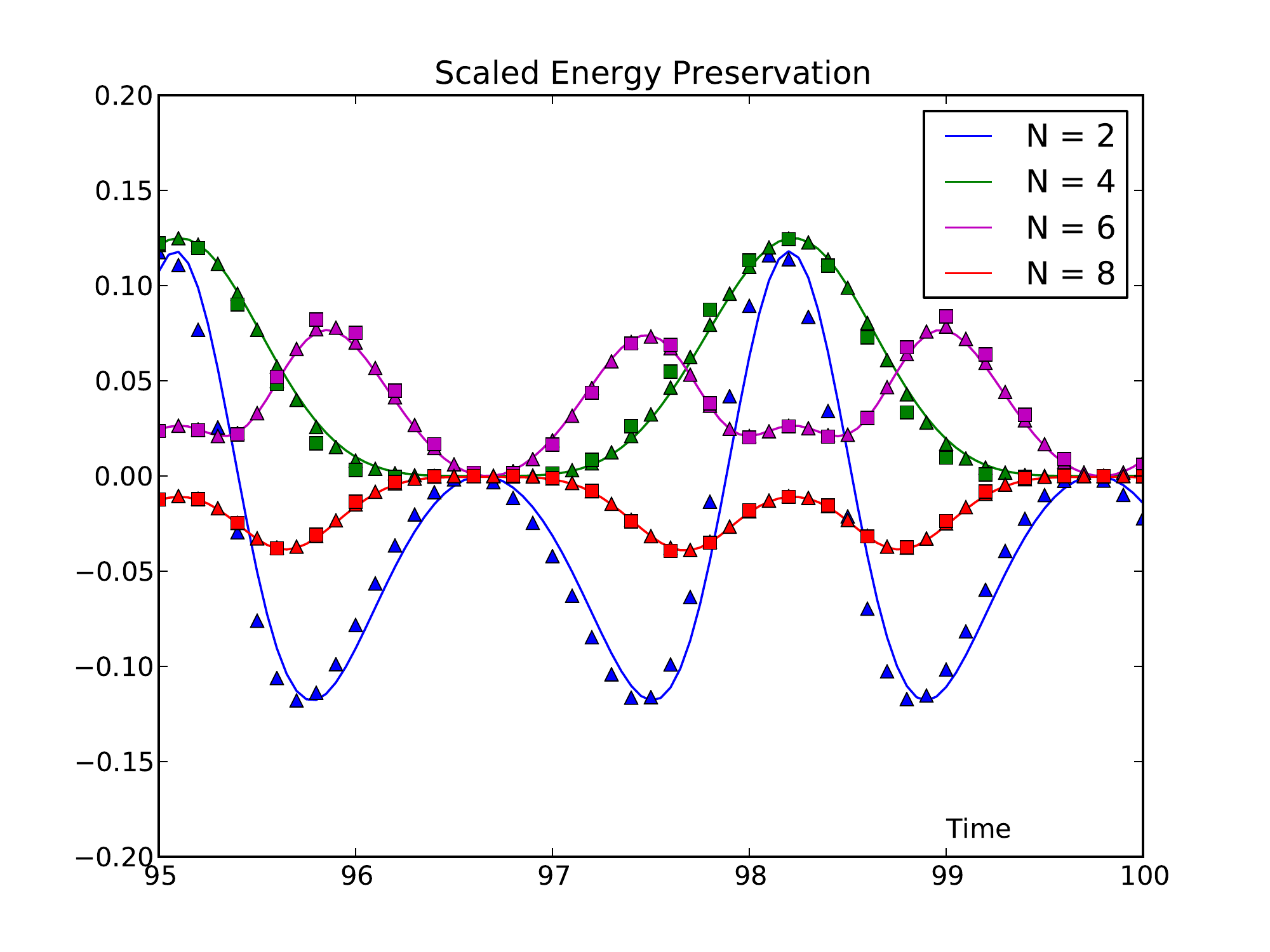}
\end{center}
\caption{This figure illustrate how well energy is conserved with the
various splitting schemes. The quantities plotted is $(H-\frac{1}{2})/\tau^{N}$
for $\tau = 0.2$ (squares), $\tau=0.1$ (triangles) and $\tau=0.05$ (lines).
Here $N=2$ for the St{\"o}rmer-Verlet scheme (blue line), $N=4$ for the
$\tau^2$-corrected generators (green line), $N=6$ for the
$\tau^4$-corrected generators (magenda line), and $N=8$  for the
$\tau^6$-corrected generators (red line). Each plotted quantity is essentially the
value of the next correction at the visited point in phase space. Since the plot is taken
over the last half of the $16^{\text{th}}$ period the figure also give some indication of how
well the exact oscillation period is reproduced by the scheme. The deviation is quite
large for the St{\"o}rmer-Verlet scheme when $\tau=0.2$; to avoid cluttering the figure
we have not included these points.
}
\label{energyPreservation}
\end{figure}
It remains to demonstrate that our algorithms can be applied to real
examples. We have considered the Hamiltonian
\begin{equation}
   H = \frac{1}{2} p^2 + \frac{1}{4}q^4,
\end{equation}
with initial condition $q(0)=0$, $p(0)=1$. The exact motion is a
nonlinear oscillation with $H$ constant equal to $\frac{1}{2}$, and period  
\begin{equation}
T = 4\,\int_0^{2^{1/4}} \frac{\sqrt{2}\, \text{d}q}{\sqrt{2-q^4}} = 2^{1/4}\,
\text{B}(\frac{1}{4},\frac{1}{2}) \approx 6.236\,339\ldots.
\end{equation}
Here $\text{B}(x,y) = \Gamma(x)\Gamma(y)/\Gamma(x+y)$ is the beta function.
In Fig.~\ref{energyPreservation} we plot the behaviour
of $\left(H-\frac{1}{2}\right)/\tau^{2+n}$ during the last half of the $16^{\text{th}}$ oscillation,
for various values of $\tau$ and corrected generators up to order $\tau^6$ (corresponding to $n=6$).

\begin{figure}
\includegraphics[clip, trim = 8ex 5.5ex 7.5ex 5ex,width=0.495\textwidth]{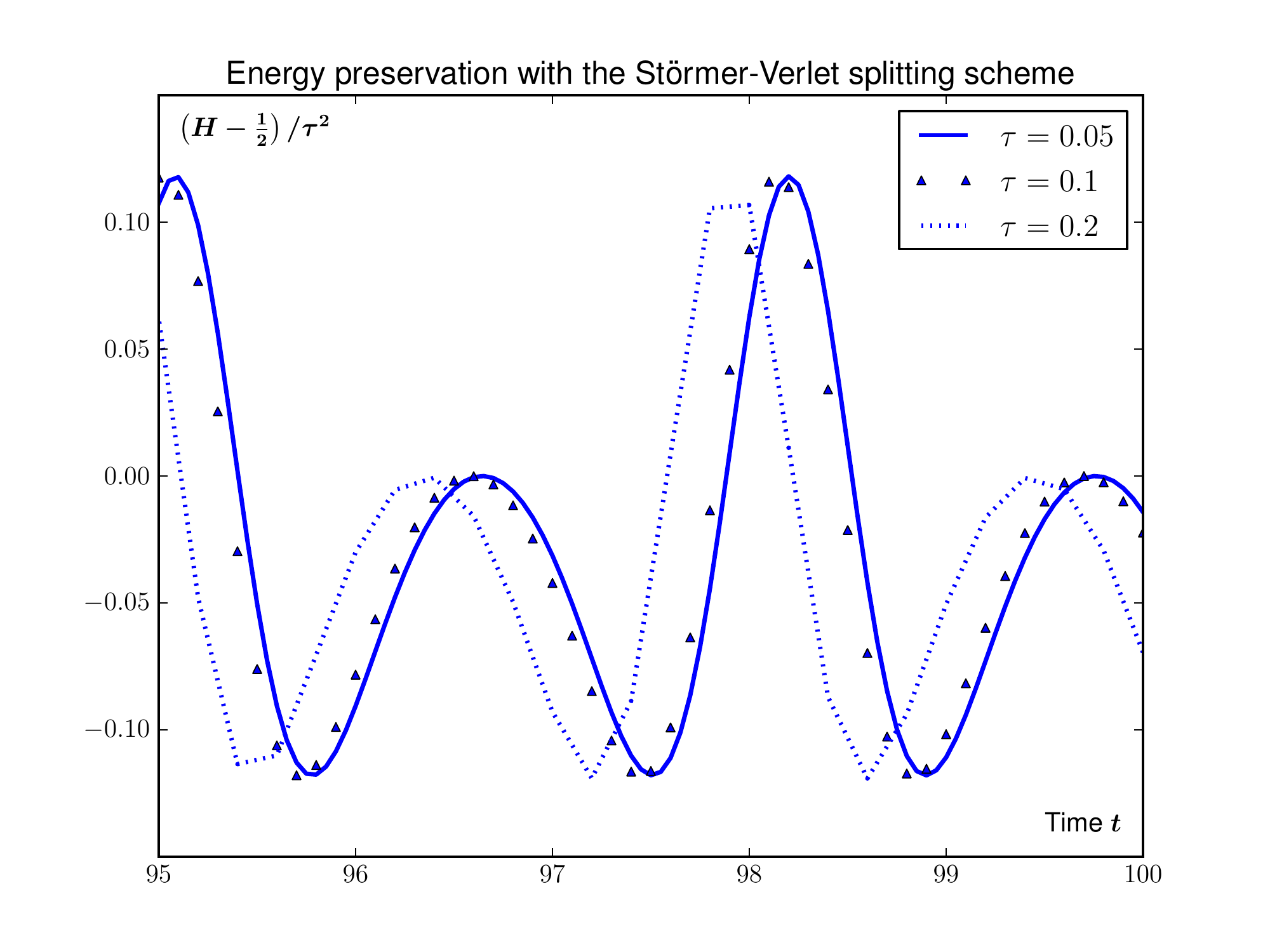}
\includegraphics[clip, trim = 8ex 5.5ex 6.5ex 5ex,width=0.495\textwidth]{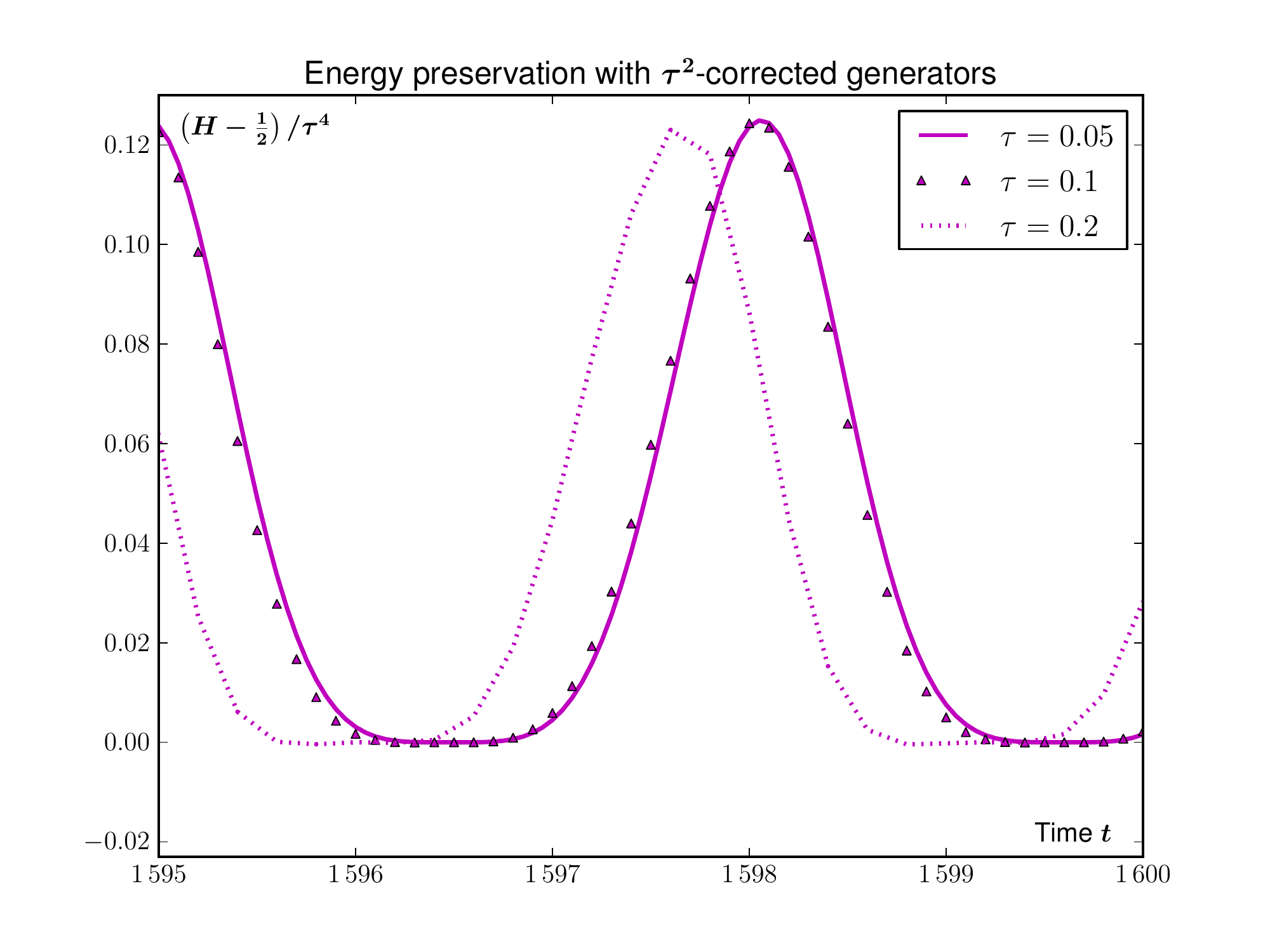}\\[1ex]
\includegraphics[clip, trim = 8ex 5.5ex 7.5ex 5ex,width=0.495\textwidth]{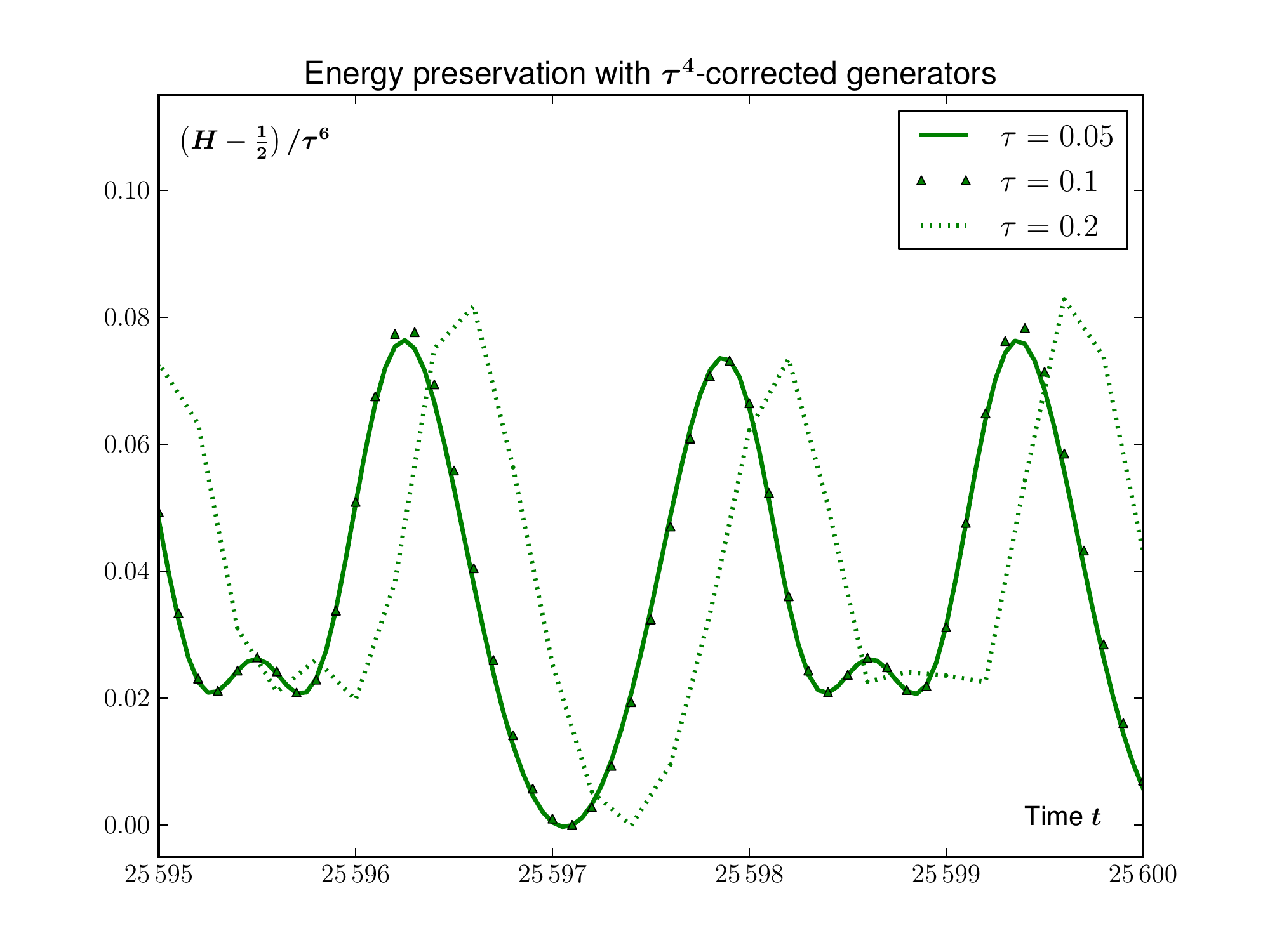}
\includegraphics[clip, trim = 8ex 5.5ex 6.5ex 5ex,width=0.495\textwidth]{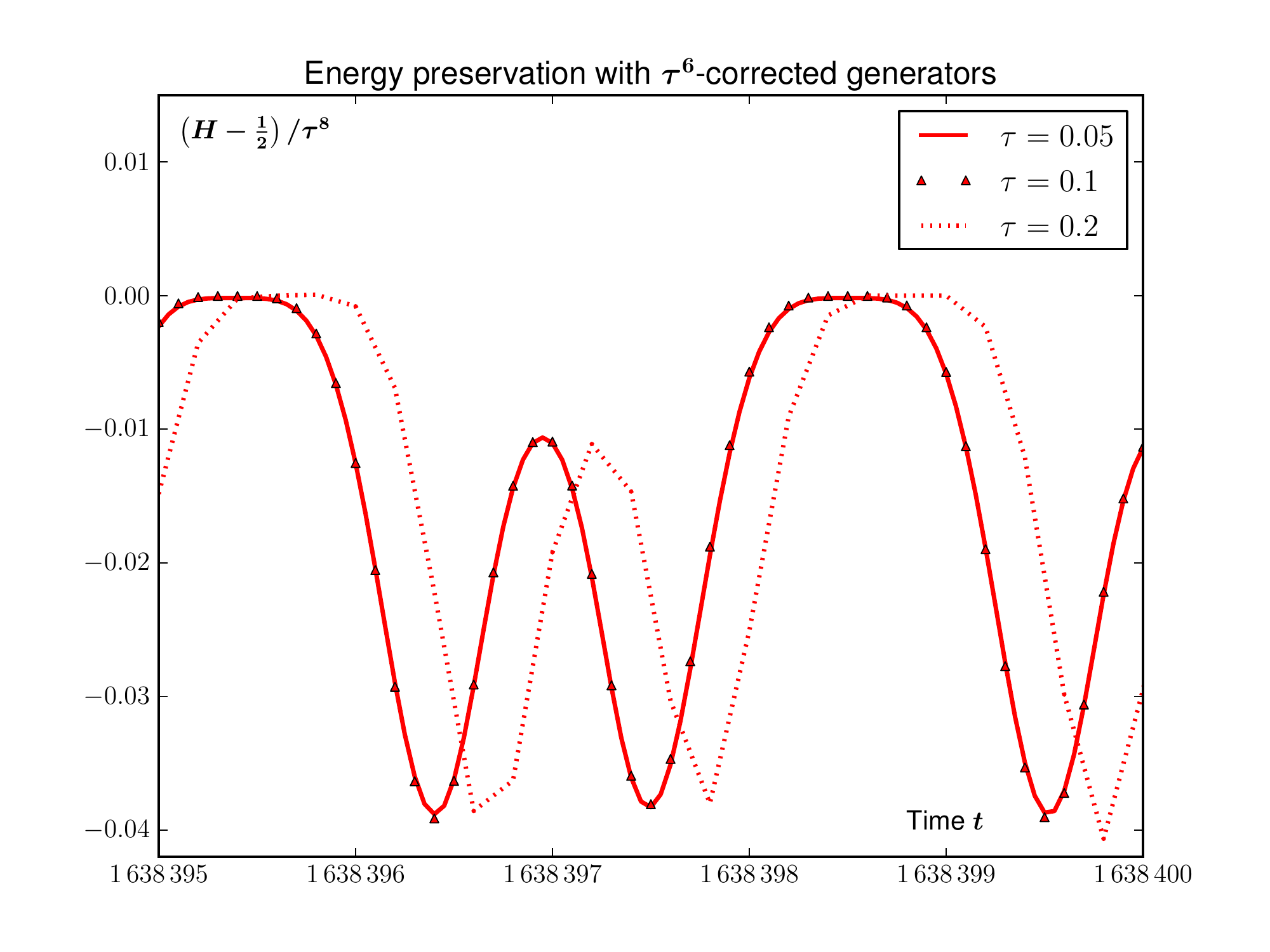}
\caption{These figures illustrate the long time behaviour through the
  last half of the $16^{\text{th}}$ period for the
  St{\"o}rmer-Verlet, first half of  $257^{\text{th}}$ period for
  $\tau^2$-corrected, last half of the $4\,104^{\text{th}}$ period for
  $\tau^4$-corrected and last half of the $262\,718^{\text{th}}$ period for
  $\tau^6$-corrected schemes. Different timesteps $\tau$ have an
  effect on the period of oscillation, but the preservation of energy remains stable for a very long time.}
\end{figure}

\subsection{Fermi-Pasta-Ulam-Tsingou type problems}

\begin{wrapfigure}{r}{0.5\textwidth}
  \vspace{-30pt}
  \begin{center}
    \includegraphics[width=0.48\textwidth]{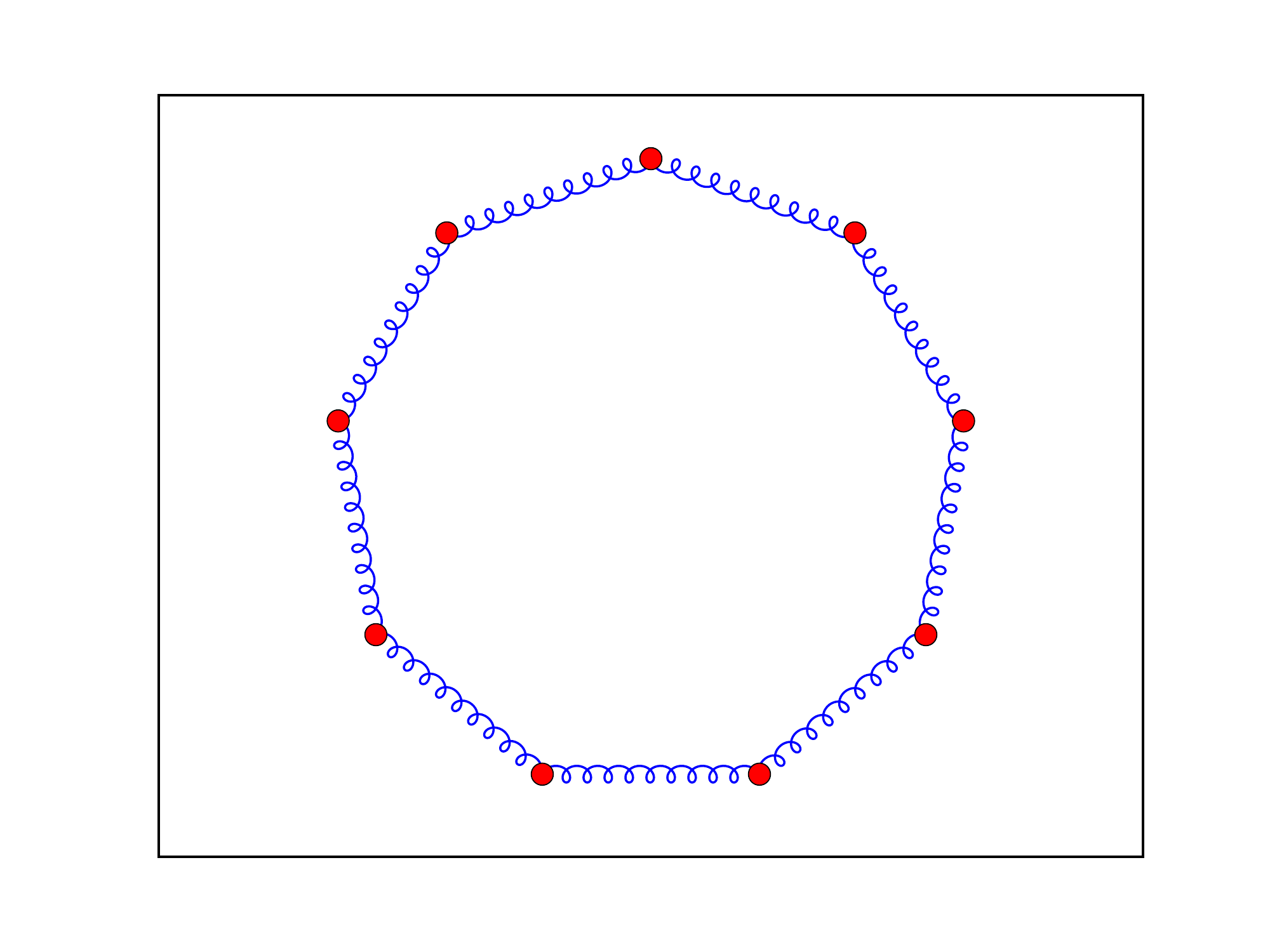}
  \end{center}
  \vspace{-30pt}
  \caption{Chain of neighbouring particles} \label{chain}
  \vspace{-10pt}
\end{wrapfigure}

We also want to demonstrate that our algorithms can be applied to
systems with many degrees of freedom, like lattice models with short
range interactions.
Here we will consider a one-dimensional closed chain of $d$
particles, as illustrated in Fig.~\ref{chain}, interacting
with its nearest neighbours through a potential $U$, and possibly with a local
substrate through a potential $V$. The latter will confine the
$n^{\text{th}}$ particle to the vicinity of a position ${\cal R}_n = n L/d$, where
$L$ is the circumference of the chain.

This class of models include the Fermi-Pasta-Ulam-Tsingou (FPU) problem
introduced in $1953$ by Fermi~\emph{et. al.}~\cite{FermiPastaUlam}
for investigating equipartition of energy among the degrees
of freedom in the system. Much research in different fields of mathematics and
physics have been devoted towards understanding the highly unexpected
dynamical behavior of this system. A review of the last $50$ years of
comprehensive study has been given by Berman and Izrailev \cite {BermanIzrailev}. 

In a recent paper Hairer and Lubich \cite{HairerLubich} presented
an investigation of the FPU problem using a modulated Fourier expansions on chains
with a large number of particles, and in \cite{McLachlanNeale} McLachlan and
Neal have made a comparison of various integrators applied to the FPU problem.
A good analysis, using the Baker-Cambell-Hausdorff formula,
of the St{\"o}rmer-Verlet method applied this problem 
has been given by Benettin and Ponno \cite{BenettinPonno}.
Application of numerical methods to this problem is also discussed by
Palearis and Penati \cite{PalearisPenati}.

Here we will demonstrate that our integrators can
be implemented and applied in practise to these type of models.
There is of course a computational cost per timestep by going
to a higher order method, but asymptotically
that cost grows linearely with the size
of the system, provided interactions are of short range.
There is also a cost in complexity 
of code implementation, which we have solved by
writing a program for automatic generation of
the numerical code \cite{AsifAnneKareIII}.
 
Let $q_m(t) = {r}_{m}(t) -{\cal R}_m$, where $r_{m}(t)$ the position of 
the $n^{\text{th}}$ particle, and consider the system described by the Hamiltonian
\begin{align}
   H({\bf q},{\bf p})=\frac{1}{2} \sum_{m=0 }^{d-1} p^2_m+  
   \sum_{m=0}^{d-1} V(q_m) +\sum_{m=0}^{d-2} U(s_m) ,\label{FPU}
\end{align}
where $d$ is the number of particles, and $s_m = q_{m+1}-q_{m}$.
A class of model which includes both the linear chain and the FPU model
can be obtained by choosing 
\begin{align*}
V(q) = \frac{1}{2}\omega^2 q^2, \quad U(s) = \frac{1}{2} s^2 + \frac{1}{3}\alpha s^3 + \frac{1}{4} \beta s^4,
\end{align*}
where $U$ is describing the interactions between particles.
This model has been referred to as the FPU 
$\alpha + \beta$ model (with $\omega^2=0$).
%
%
In this example we have used $\alpha =0$ and $\beta = 1$.
We have tested the methods with respect to (i) energy conservation,
(ii) deviation of the generated solution from the exact
solution\footnote{Actually a numerical solution of the same system computed
to very high precision.}, and (iii) efficiency of the methods with respect to CPU time.

\begin{figure}
\includegraphics[clip, trim = 7.5ex 5.5ex 10.5ex 1.5ex,width=0.99\textwidth]{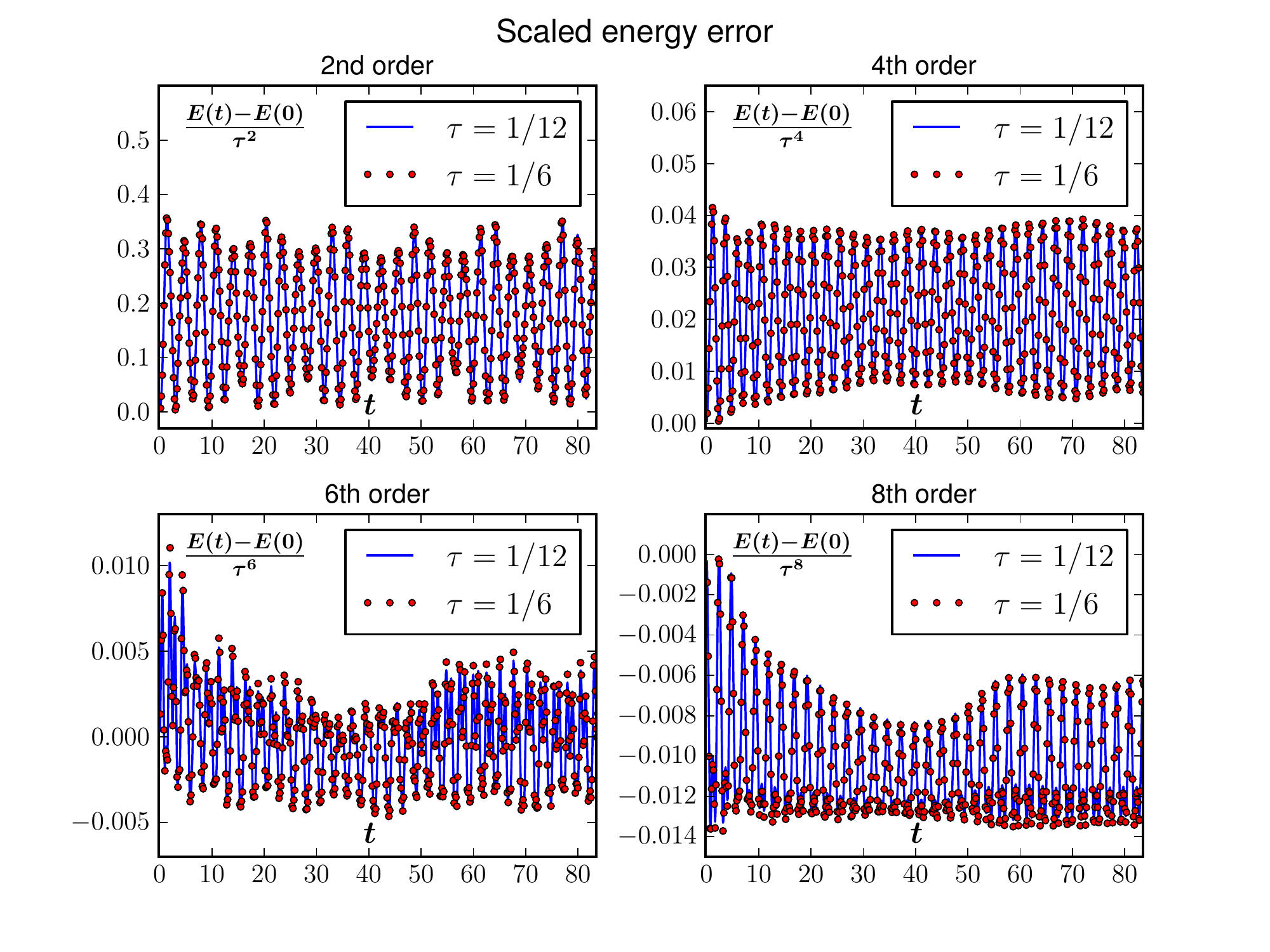}
\caption{Scaled energy error for higher order methods of the FPU} \label{energyError}

\end{figure}
\begin{figure}
\includegraphics[clip, trim = 7.5ex 5.5ex 10ex 1.5ex, width=0.99\textwidth]{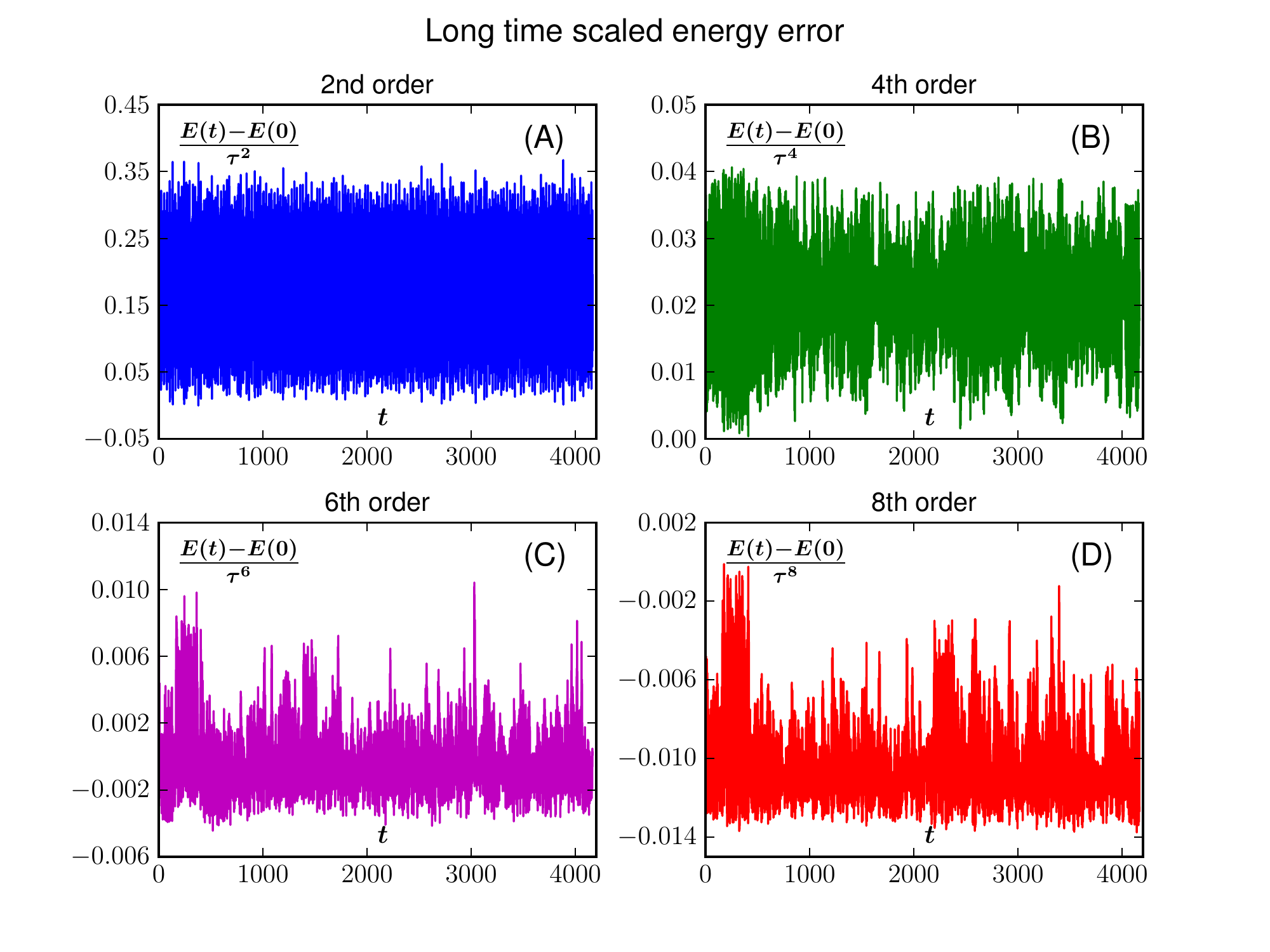}
\caption{Long time scaled energy error for a Fermi-Pasta-Ulam-Tsingou type problem, computed with
the St{\"o}rmer-Verlet (A) and higher order corrected integrators (B--D) with
timestep $\tau=\frac{1}{12}$ and initial energy $E(0)=1.425$.} \label{LongtimeenergyError}
\end{figure}

In Fig.~\ref{energyError} we show the scaled energy error on FPU for
different choices of $\tau$ with all four methods. For these experiments
we consider $9$ particles with initial energy $E(0)=1.425$.
As can be seen the energy conserved very well for all methods, with
the error scaling like $\tau^N$ for a method of order $N$. 
As demonstrated by the long time behaviour in Fig.~\ref{LongtimeenergyError}
the energy error does not increase noticably with time.

Another quantity of interest in a system with many degrees of freedom
is the \emph{global error}, i.e.~a measure how much the numerical
solution deviates from the exact solution. Here an exact solution
is not available. Instead we have generated a very accurate solution
by use of our eight order method with timestep
$\tau=5\cdot 10^{-4}$, calculated with multiprecision (50 decimal digits)
floating point accuracy. This is for practical purposes as good as
an exact result, and we will refer to it as such.

\begin{center}
\begin{figure}
\centering
\includegraphics[clip, trim = 3.5ex 3.5ex 9ex   3.5ex, width=0.5\textwidth]{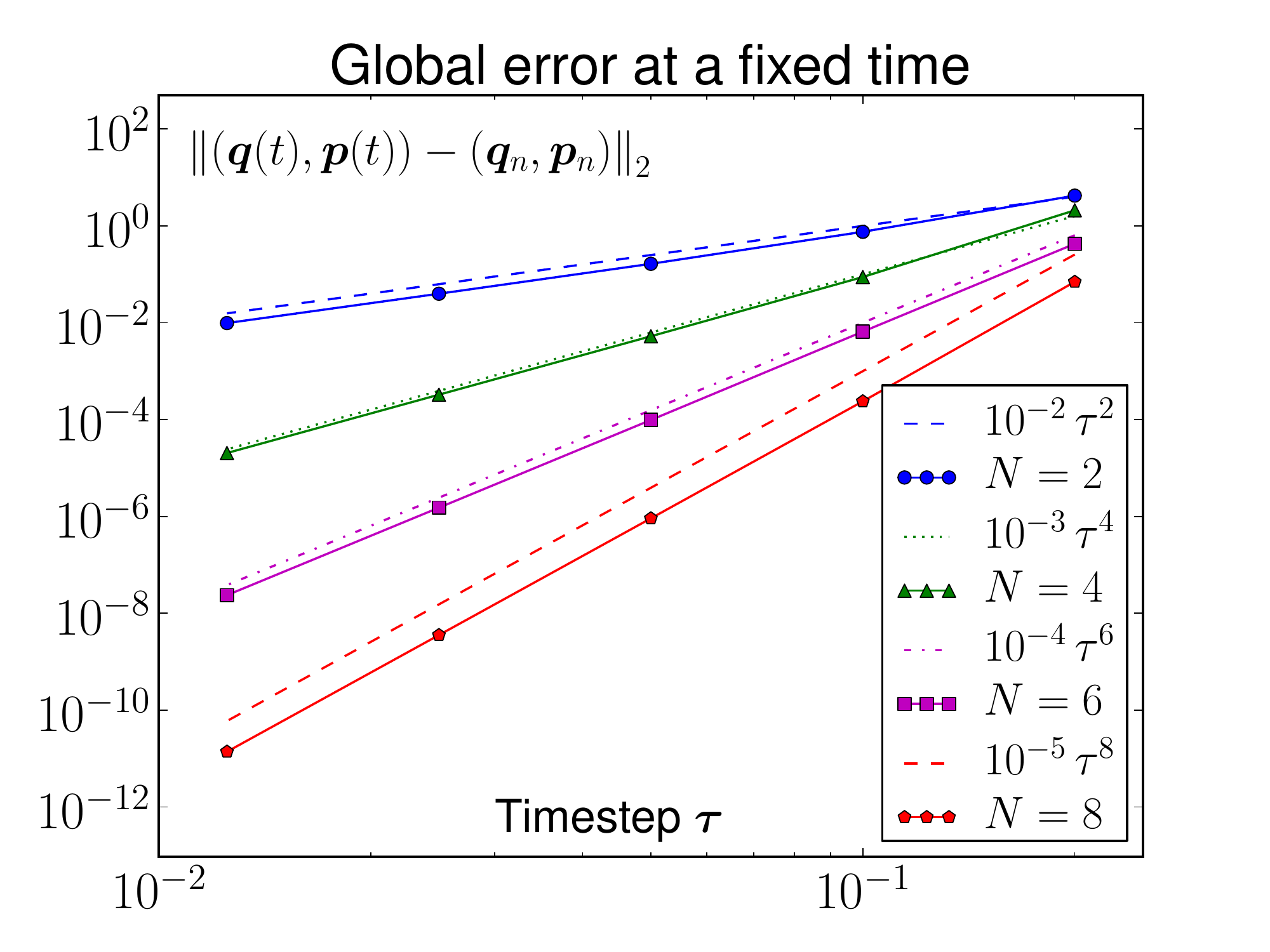}\hspace{0.0em}
\includegraphics[clip, trim = 3.5ex 3.5ex 8ex 3.5ex,width=0.5\textwidth]{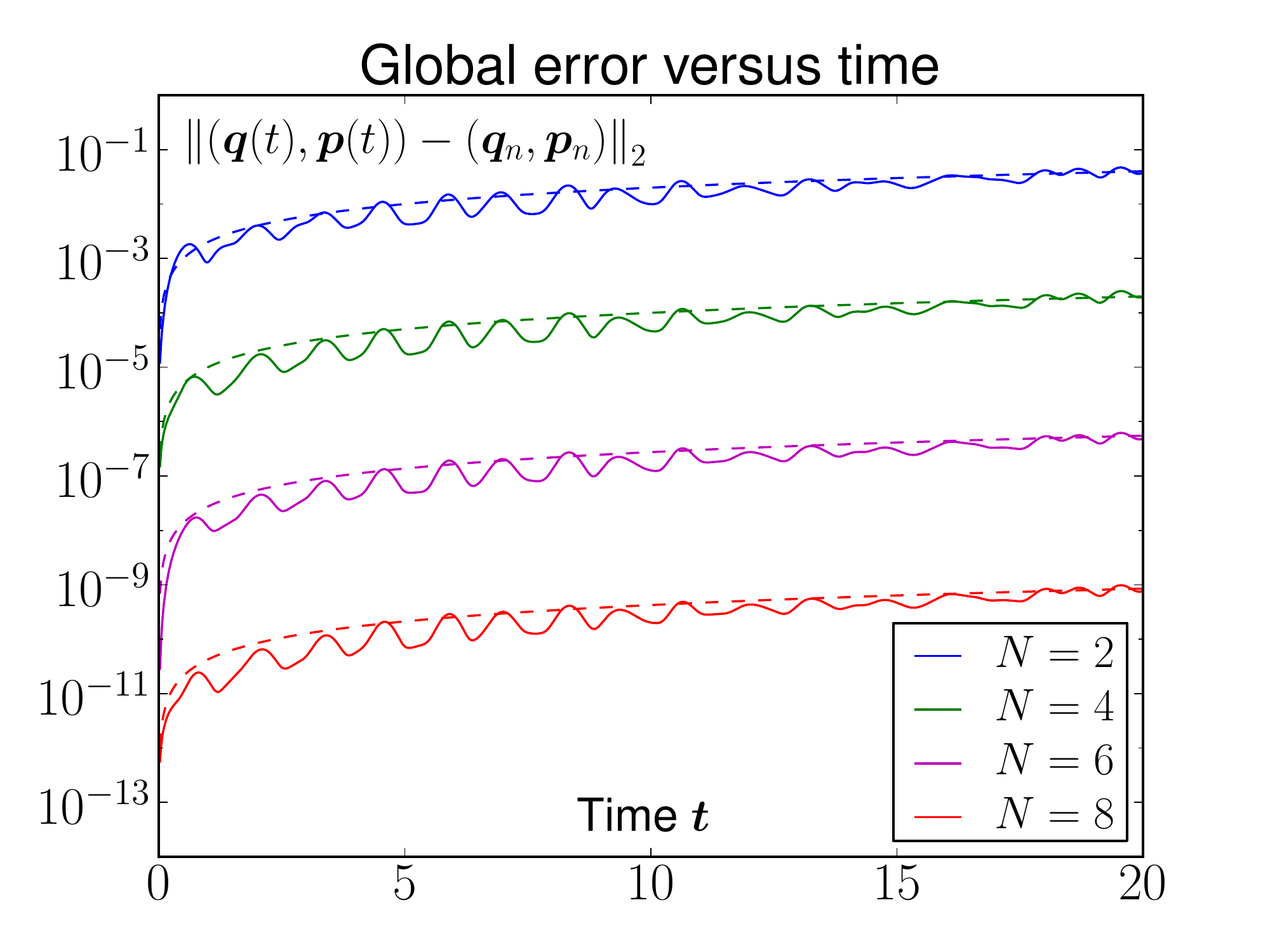}
\caption{The left frame shows how the global error (here measured at time $t=10$) depends
on the timestep $\tau$ and the order $N$ of the integration scheme. As expected this error varies
like $\tau^N$ (as long as it is small).
The right frame shows how the global error grows with time,
here for a timestep $\tau=\frac{1}{40}$. The dashed lines are eyeball fits to linear error growth,
cf.~$\varepsilon(t) \equiv \left\| (\bm{q}(t), \bm{p}(t)) - (\bm{q}_n, \bm{p}_n) \right\|_2
\sim C\, t\, \tau^{N}$. These results are for a lattice of $d=9$ particles.
} \label{globalError}
\end{figure} 
\end{center}

We have investigated several measures of deviation;
they all give qualitatively the same results. Here
we will only discuss the quantity
\begin{align}
    \varepsilon(t) &\equiv \left\| (\bm{q}(t), \bm{p}(t)) - (\bm{q}_n, \bm{p}_n) \right\|_2\nonumber\\
    &= \left[ \sum_{m=0}^{d-1} (q_m(t) -q_{m,n})^2 + (p_m(t) -p_{m,n})^2 \right]^{1/2},
\end{align}
where $\bm{q}_n$ ($\bm{p}_n$) denote the positions (momenta) of the numerical solution
at a timestep $n$ such that $n\tau = t$, and $\bm{q}(t)$ ($\bm{p}(t)$) denote the positions (momenta)
of the exact solution at time $t$. As shown in Fig.~\ref{globalError} the global error behaves roughly
like
\begin{equation}
   \varepsilon(t) \sim C\, t\, \tau^{N}, \label{GlobalErrorBehaviour}
\end{equation}
for relatively short times $t$. Here $C$ is a constant which depends
on the order $N$ of the method and the initial conditions. This is
in agreement with exact behaviour of integrable
systems, cf.~Theorem 3.1 in the book \cite{HLWG} by Hairer \emph{et.~al.}

\begin{center}
\begin{figure}
\centering
\includegraphics[clip, trim = 5ex 3.5ex 8ex
3.5ex,width=0.8\textwidth]{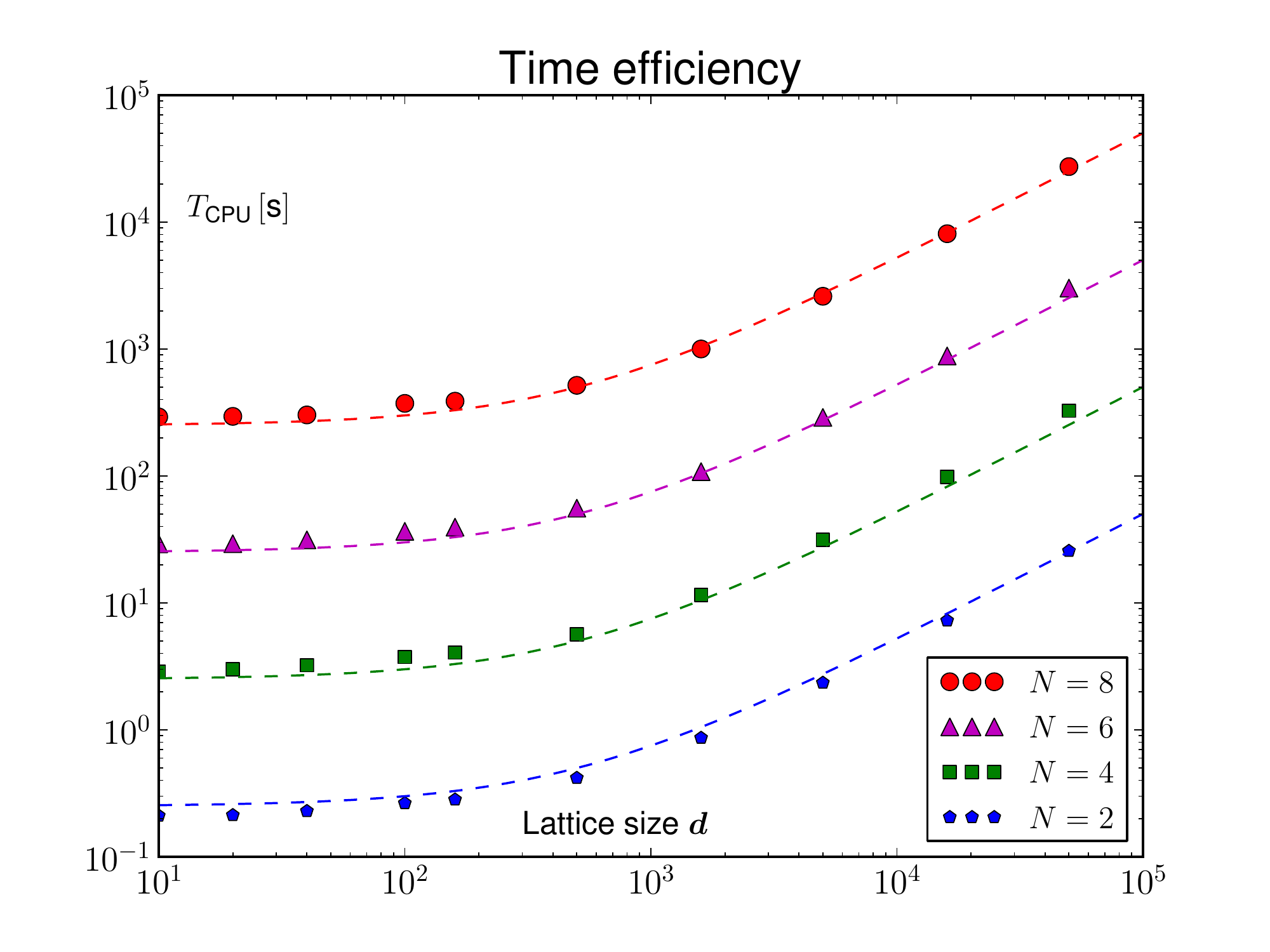}
\caption{CPU time $T_{\text{CPU}}$ used to solve a lattice of $d$ particles for $1000$ timesteps
($\tau=\frac{1}{12}$) for schemes of different orders $N$. Asymptotically, $T_{\text{CPU}}$ grows
linearly with $d$. The penalty for increasing the order $N$
by 2 is about a factor $10$ increase in $T_{\text{CPU}}$,
when $d$ and the number of timesteps is kept fix.
} \label{CPUtime}
\end{figure} 
\end{center}

To check the efficiency of our methods in practical use,
we have also measured CPU time used to integrate systems with different number $d$ of particles,
with $d$ ranging from $9$ to $50\,000$. All runs have been done on the same system,
a workstation equipped with two six-core Opteron 2431 processors, but
using code written in \texttt{NumPy}. Hence, the code
is not parallellized and run on a single core.
Some results, run with timestep $\tau = 1/12$ for all methods, 
is shown in Fig.~\ref{CPUtime}. Under these conditions we find that
the CPU time increases by a factor of about $10$ for each step in
order. From the left frame of Fig.~\ref{globalError} we see that
this step also increases the accuracy with a factor of
about $10^{-1}\,\tau^2$ (for $d=9$ particles). If we want a
prescribed accuracy $10^{-P}$ for the global error $\varepsilon(t)$ at time $t$
we may choose to use lower order method with a small timestep
(which requires many steps $n$), or a higher order method with
fewer, but more time-consuming steps. Which choice is best?
For the parameters displayed in Fig.~\ref{globalError} we estimate
the condition
\begin{equation}
     \varepsilon(t) \approx 10^{-2-N/2}\,t\,\tau^N \approx 10^{-P}.
\end{equation}
I.e., we must choose a timestep such that
\begin{equation}
     \tau^N \approx \frac{1}{t}\times10^{2+N/2 -P},
\end{equation}
which requires
\begin{equation}
     n \approx \frac{t}{\tau} \approx \frac{1}{\sqrt{10}}\times t^{(1+1/N)}\,10^{(P-2)/N}
\end{equation}
steps, where each steps requires a CPU time
$t_{\text{step}} \approx 10^N\, t_0$ for some constant $t_0$ which depends on the computer being used.
Hence, we should choose $N$ to minimize
\begin{equation}
     T_{\text{CPU}} = n\,t_{\text{step}} \approx \, \frac{t_0\,t}{\sqrt{10}}\,\times \,10^{N/2 +(P+\log_{10} t-2)/N}.
\end{equation}
Treating $N$ as a continuous varible gives the optimal value
\begin{equation}
     N_{\text{opt}} \approx \sqrt{2\,\left(P + \log_{10} t -2\right)}.
     \label{OptimalOrder}
\end{equation}

\section{Concluding remarks}\label{concludingRemarks}

In this paper we have shown that it is possible to systematically
extend the standard St{\"o}rmer-Verlet symplectic integration scheme
to higher orders of accuracy, and that the higher order schemes can
be applied in practise to physical systems of interest, including FPU-type 
lattice problems with many particles (with nearest-neighbour interactions).
As illustrated by equation~(\ref{OptimalOrder}),
it is advantageous to use a higher order method when one wants a solution of high
precision $P$, and also if one wants a solution of moderate accuracy but
over a long time interval.

As demonstrated, the theoretical algorithms have been implemented and tested.
One rapidly discovers that it is a nightmare to do a correct implementation by hand.
The general compact form of these schemes usully expand to very long
expressions, which are laborious and error-prone to handle manually.
We have therefore developed a set of computer routines
which automatically generate the basic numerical integrators
for a complete timestep of each specific model.

For the cases we have investigated these integrators perform according to
expectations, sometimes even better than expected.

\section*{Acknowledgements}

We thank professor Ernst Hairer for encouraging remarks
and helpful pointers to the literature.

\newpage

\end{document}